\documentclass{article}
\usepackage{fancyhdr}
\RequirePackage[OT1]{fontenc}
\usepackage{amsthm,amsmath,natbib,bm}
\RequirePackage[colorlinks,citecolor=blue,urlcolor=blue]{hyperref}
\usepackage{epsfig,multicol,multirow}
\usepackage{geometry}
\usepackage{graphicx}
\usepackage{float}
\usepackage{placeins}
\usepackage{booktabs}
\usepackage[table]{xcolor}
\usepackage{adjustbox}
\usepackage{econometrics}
\def\Tf{\mathcal{T}}

\def\MADf{\mathcal{MAD}}

\def\IF{\hbox{IF}}
\def\R{\mathcal{R}}
\def\D{\mathcal{D}}

\newcommand{\PIF}{\text{PIF}}

\def\mad{\hbox{MAD}}

\def\asd{\hbox{ASD}}

\def\asv{\hbox{ASV}}

\def\wasd{\widehat{\asd}}

\newcommand{\asim}{\stackrel{\text{\tiny approx.}}{\sim}}

\numberwithin{equation}{section}
\theoremstyle{plain}
\newtheorem{theorem}{Theorem}[section]

\fancypagestyle{firstpage}{
    \fancyhf{} 
    \fancyhead[C]{This is a earlier version of a paper that has now been published [Communications in Statistics - Simulation and Computation, 1–10]. \\ Please cite as: Arachchige, C. N. P. G., \& Prendergast, L. A. (2024). Confidence intervals for median absolute deviations. Communications in Statistics - Simulation and Computation, 1–10. https://doi.org/10.1080/03610918.2024.2376198}
}

\begin{document}

\title{Confidence intervals for median absolute deviations}
\author{Chandima N. P. G. Arachchige \\ Department of Mathematics and Statistics, La Trobe University \\ \href{mailto:
18201070@students.latrobe.edu.au
}{
18201070@students.latrobe.edu.au
}\\ \\ Luke A. Prendergast\\ Department of Mathematics and Statistics, La Trobe University, \\ Australia\\ \href{mailto:luke.prendergast@latrobe.edu.au}{luke.prendergast@latrobe.edu.au}}
\date{}

\maketitle
\thispagestyle{firstpage}

\begin{abstract}
The median absolute deviation (MAD) is a robust measure of scale that is simple to implement, easy to interpret and often reported as a measure of dispersion when data is skewed. However, confidence intervals either for a single MAD or the differences/ratio of two MADs are lacking, meaning that interpretations may be subjective.
Motivated by this, we introduce interval estimators of the MAD to make reliable inferences for dispersion for based on MADs. Our simulation results show that the coverages of the intervals are very close to the nominal coverage for a variety of distributions, including those that are heavily skewed.  We also present examples that clearly highlight the advantages of using these intervals instead of the heavily criticised intervals derived from sample variances when the data is skewed. 
\end{abstract}


\section{Introduction} \label{sec:Intro}
The median absolute deviation  is a robust measure of dispersion \citep[MAD, see e.g.][]{hamp-1974,hampel1986robust}. Defined as the median of the absolute residuals from the median, the MAD is a suitable scale measure to accompany the median. \cite{hamp-1974} referred to the $\mad$ as the ``median deviation" and it had first received attention even as early as \cite{gauss1816demonstratio}, and was later rediscovered by \cite{hampel1968contribution}. The MAD is the \textit{most robust} estimator of scale  as measured by robustness measures such as the break-down point and gross error sensitivity \citep{hamp-1974}. The breakdown point of an estimator is the proportion of contamination that the estimator can handle before providing  unreliable results and for the MAD this is equal to 1/2 (the maximum).  The $\mad$ has what is known as a bounded influence function so that the amount of influence any observational type can exert on the estimator is limited.  Therefore, the MAD is very suited to skewed data or data including outliers where, in both cases, the use of the sample standard deviation can be misleading and associated tests unreliable (CITE). 

Despite the fact that the MAD is often reported, it is done so as a point estimate with little or no consideration for uncertainty in estimation.  

\cite{arachchige2019interval} showed that excellent coverages for interval estimators of ratios of interquantile ranges can be achieved.  This makes these intervals more suitable than those for ratio of variances when normality cannot be assumed.  Then, \cite{arachchige2019CV} considered interval estimators for robust versions of the coefficient of variation, one of which uses the MAD in place of the standard deviation (and the median to replace the mean). Motivated by these good coverage properties, we consider interval estimators for the MAD and  for ratios and differences of independent MADs as robust alternatives to intervals based on sample variances. To the best of our knowledge, and not to confuse the MAD with the \textit{mean} absolute deviation for which interval estimators with good coverage have been introduced by \cite{bonett2003confidence}, no one has introduced these interval estimators for the MAD.  The very good coverage properties, that we will highlight later, ensure inferences about dispersion based on the MAD are possible. Based on a preprint of our work in \cite{arachchige2019confidence}, our confidence interval for MAD has been included in the \textit{DeskTools} 
\citep{signorell2019desctools} package in R.

In Section 2 we provide some necessary notations before considering influence functions for ratios of MADs.  In Section 3 we consider confidence intervals for MADs, differences of MADs and ratios of MADs with coverage properties explored via simulations in Section 4.  Examples are also considered in Section 4 and we conclude in Section 5.
\section{Notations and influence functions}\label{sec:MAD}
Let $X$ denote a random variable and $F$ its distribution function.  Then \cite{hamp-1974} defined the median absolute deviation (MAD) as
\begin{equation}
 \mad(X) =\text{med}\mid{X-M}\mid~, 
 \label{eq:MAD}
\end{equation}
where `med' denotes the median and $M=\text{med}(X)=F^{-1}(0.5)$ is the population median.  Let $X_1,\ldots,X_n$ denote a random sample of $n$ observations.  Then the MAD estimate is simply the median of the absolute residuals from the sample median.  That is, for $m$ denoting the sample median, $\widehat{\text{MAD}}$ is the sample median of the $|X_1-m|,\ldots,|X_n-m|$.  While inference, for a single MAD may be of interest, it is often that case that comparison of dispersion measures, such as the MAD, is needed to compare two populations.

Consider two independent random variables $X\sim F_1$ and $Y\sim F_2$ and let us consider $\mad(X)$ and $\mad(Y)$. Then, the 
population squared ratio of MADs, which we denote as $R_M$, and associated estimator can be define as
\begin{equation}
R_M = \left[\frac{\mad(X)}{\mad(Y)}\right]^2\;\;\text{and}\;\; \widehat{R}_M = \left[\frac{\widehat{\mad(X)}}{\widehat{\mad(Y)}}\right]^2.
\label{eq:R_m}
\end{equation}

Here we have suggested the squared ratio of MADs since it is the analogue to the ratio of variances and, in fact, equal to ratio of variances for some distributions (e.g. normal). However, the ratio of MADs may also be used.  Another possibility is the difference of MADs, $D_M$, where
\begin{equation}
 D_M = \mad(X) - \mad(Y)  \;\;\text{and}\;\;  \widehat{D}_M =  \widehat{\mad(X)} - \widehat{\mad(Y)}
\label{eq:D_m}~. 
\end{equation}
\subsection{Influence function and partial influence functions} \label{sec:IFandPIF}
Define the contamination distribution to be $F_{\epsilon}=(1-\epsilon)F+\epsilon \Delta_{x}$, where $\epsilon \in [0,1]$ is the proportion of contamination and $\Delta_{x}$ has all of its mass at the contaminant $x$. Consider an estimator functional $\Tf$ such that $\Tf(F)=\theta$ and $\Tf(F_n)=\widehat{\theta}$ where $F_n$ denotes the empirical distribution function for sample of $n$ observations. The relative influence on $\Tf(F)$ of $\epsilon$ proportion of contaminated observations at $x$ is given by,  $[\Tf(F_\epsilon)-\Tf(F)]/\epsilon$, where $\Tf(F_{\epsilon})=(1-\epsilon)\Tf(F)+\epsilon \Delta_{x}$. Then, the influence function \citep[IF][]{hamp-1974} is defined as,
 $$\IF (x;\Tf,F)=\lim_{\epsilon\downarrow 0} \frac{\Tf(F_\epsilon) - \Tf(F)}{\epsilon}\equiv \frac {\partial }{\partial \epsilon }\Tf(F_\epsilon)\Big|_{\epsilon = 0}.$$

When more than one population exists, the IF is determined by contaminating one population while the other population remains uncontaminated. \cite{pires2002partial} defines this notion as ``partial IFs" (PIFs) and in our context with two populations we have two PIFs. The first PIF of the estimator with functional $\Tf$ at ($F_1,\ F_2$) is
\begin{align}
\PIF_1 (x;\Tf,F_1,F_2) = \lim_{\epsilon \to   0}\left[\frac{\Tf[(1-\epsilon)F_1+\epsilon\Delta_{x_0},F_2]-\Tf(F_1,F_2)}{\epsilon}\right]  
\label{eq:PIF}
\end{align}
and with $\PIF_2 (x;\Tf,F_1,F_2)$ defined similarly.

Now, consider the functional for the standardized MAD denoted by $\MADf$ so that $\MADf(F)=\mad_X$. \cite{hamp-1974} gives the influence function for the $\mad$ when $F$ is the normal distribution and further details can be found on  page 107 of \cite{hampel1986robust}.  Let $f=F'$ denote the density function then, assuming $f(M)$ and $2[f(M+\text{MAD}_X)+f(M-\text{MAD}_X)]$ are nonzero, a general form of the IF for the MAD exists; e.g. see page 137 of \cite{huber1981robust} or page 16 of \cite{andersen2008modern}.  This is given as
\begin{equation}
\IF (x;\,\MADf ,F) = \displaystyle\frac{\left[\text{sign}(x-M)-\text{MAD}_X\right]-\displaystyle\frac{f(M+\text{MAD}_X)-f(M-\text{MAD}_X)}{f(M)}\text{sign}(x-M)}{2[f(M+\text{MAD}_X)+f(M-\text{MAD}_X)]}. \label{eq:IFmad3} 
\end{equation}

\subsubsection{Partial influence functions of the difference and squared ratio of MADs}\label{sec:PIF_Rm}
 Let $\D_M$ be the functional for the difference of $\mad$s so that,
\begin{equation*}
\D_M(F_1,F_2)= \MADf(F_1)-\MADf(F_2) ~    
\end{equation*}
then the PIFs are $\PIF_1 (x;\D_M,F_1,F_2) =  \IF (x;\MADf,F_1) $ and $\PIF_2 (x;\D_M,F_1,F_2) = -\IF (x;\MADf,F_2)$.  These are trivial and previous studies on robustness of the MAD may be considered for this context.  We therefore do not explore the difference PIFs further.

Let $\R_M$ be the functional for the squared ratio of $\mad$s so that,
\begin{equation*}
\R_M(F_1,F_2)=\left[\frac{\MADf(F_1)}{\MADf(F_2)}\right]^2 ~.    
\end{equation*}
Then the PIFs for the squared ratio of MADs are given below.
\begin{theorem}
For $\PIF (x;T,F_1,F_2)$ as defined in \eqref{eq:PIF}, the PIFs of $\R_M$ are
\begin{align*}
\PIF_1 (x;\R_M,F_1,F_2) &=  \frac{2\,\R_M(F_1,F_2)}{\MADf(F_1)}\IF (x;\MADf,F_1) ,\nonumber\\
 \PIF_2 (x;\R_M,F_1,F_2) &= -\frac{2\,\R_M(F_1,F_2)}{\MADf(F_2)}\IF (x;\MADf,F_2).
 \end{align*}
 \label{th:PIF_Rm} 
\end{theorem}
\noindent The proof of Theorem \ref{th:PIF_Rm} is in Appendix \ref{app:proof_of_PIF_Rm} and we consider some examples of the first PIF next.
 
\subsubsection{Partial influence functions comparison}\label{sec:PIF_Comp}
 \begin{figure}[!htb]
 \centering
  \includegraphics[width=\linewidth]{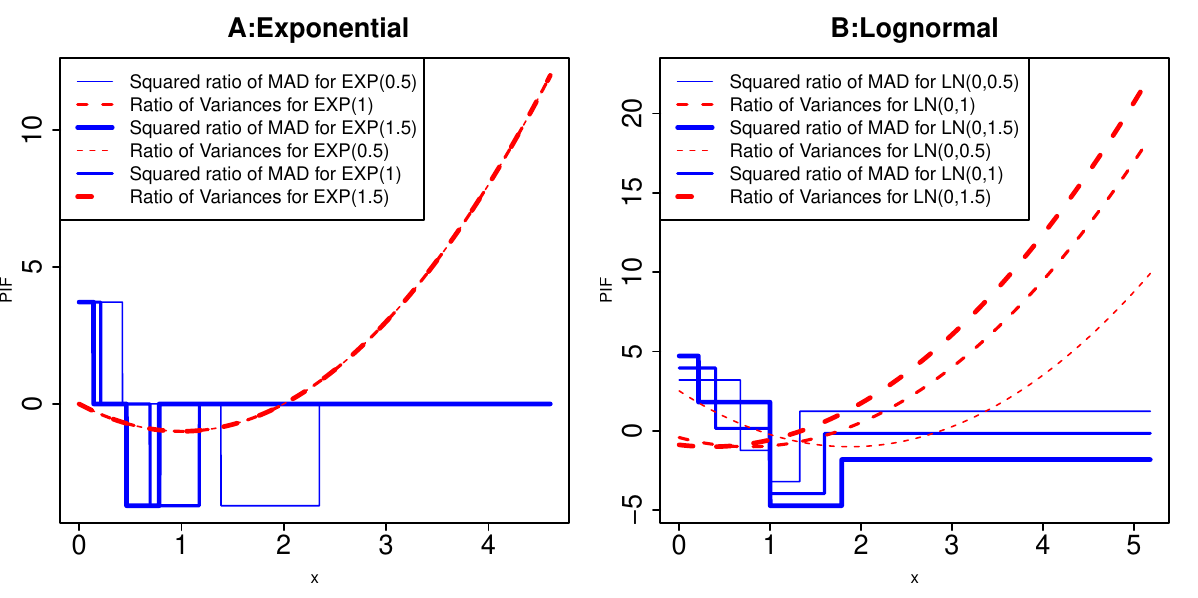}
  \caption{PIF$_1$ comparisons for (A) two exponential populations both with rates 0.5,  1 and 1.5 and (B) two log-normal populations both with $\mu$=0 and $\sigma$=0.5,1,1.5.}
  \label{fig:PIF_Com}
\end{figure}
Figure \ref{fig:PIF_Com} depicts the $\PIF$s of the first population for the squared ratio of $\mad$s and the ratio of variances  \citep[see][for these]{arachchige2019interval}.  In Plot A we consider the ratio of variances and squared ratio of MADs for two exponential distributions, both with rates equal to 0.5, or 1 or 1.5. Similarly, in Plot B we do this for two log normal distributions both with $\mu$=0 and $\sigma=0.5$ or 1 or 1.5. Since the numerator and denominator distributions are the same, both are estimators of one and therefore the PIFs are comparable.  As expected, the $\PIF$s of the ratio of variances is unbounded indicating that outliers can exert large influence on the estimator. The $\PIF$s of the squared ratio of MAD is bounded and the influence of any large outliers is limited, and far less than for the ratio of variances. For the exponential distribution, the $\PIF$s of ratio of variances do not depend on the rate parameter. However, for the log-normal distribution the PIF for the ratio of variances increases quickly with increasing $\sigma$. 

\section{Asymptotic confidence intervals} \label{sec:CIs}

In their discussion of intervals for the mean absolute deviations, \cite{bonett2003confidence} provide suggestions for median absolution deviations from a fixed point, $h$.  They suggest using intervals for the median and where the data used is the transformed $|X_i-h|$s. When $h$ is the population median, i.e. $h=M$, and this median is known, simulations (Table \ref{tab:CI_MADfromTarget1}, Table \ref{tab:CI_ratioMADfromTarget} ) result in good coverage that is close to nominal. However, when $M$ is not known and needs to be estimated, this approach typically results in coverage that is too low (e.g. less than 0.8 for a nominal 0.95). In this section we therefore provide confidence intervals that have good coverage properties, as shown by our simulations that follow.

Asymptotic normality and associated variance of the $\mad$ can be found in \cite{falk1997asymptotic} who provide the asymptotic joint normality between the median and $\mad$ estimators.  We again let $\text{MAD}_X=\MADf(F)$ and also let $\MADf(F_{n}) = \widehat{\text{MAD}}_X$.  Then, if $F$ is continuous near, and differentiable at, the median $M$, $M-\mad_X$ and $M + \mad_X$ with $f(M)>0$ and $ B_1 =f(M - \mad_X) + f(M + \mad_X) > 0,$ we have
\begin{equation*}
  \sqrt{n}\left(\widehat{\text{MAD}}_X-\mad_X\right) \asim  N(0, \asv)~, 
   \label{eq:ASyNormal_mMAD}
\end{equation*}
where `$\asim$' denotes `approximately distributed'.  The asymptotic variance of the MAD estimator is
\begin{align}
\asv = \asv(\MADf; F)=\frac{1}{4B_1^2}\left[1+\frac{B_2}{\left[f(M)\right]^2}\right] ~,
\label{eq:ASV_MAD}
\end{align}
where $B_1$ is given above and $B_2=B_3^2 + 4B_3f(M) \left[1-F(M+\mad_X)-F(M-\mad_X)\right]$ with $B_3 =f(M-\mad_X)-f(M+\mad_X)$.

We used the $\asv$ in \eqref{eq:ASV_MAD} and the Delta method \citep[see e.g., chapter 3 of ][]{Das-2008} to derive the asymptotic variance of the  ratios of $\mad$s. The asymptotic variance of $\sqrt{n_1 + n_2}\R_M(F_{n_1}, F_{n_2})$ is
\begin{equation}
\asv (\R_M;n_1,n_2)= 4\R^2_M(F_1,F_2)\Bigg[\frac{\asv(\MADf,F_1)}{w_1\,\MADf^2 (F_1)}+\frac{\asv(\MADf,F_2)}{w_2\,\MADf^2 (F_2)} \Bigg]~
\label{eq:ASV_Rm}  
\end{equation}
where $w_i=n_i/(n_1+n_2)$ for $i=1,2$.

Since the two populations are independent, deriving the asymptotic variance of the difference of MAD is straightforward.
\begin{equation}
\asv (\D_M;n_1, n_2)=  \asv(\MADf,F_1) + \asv(\MADf,F_2).
\label{eq:ASV_Dm}  
\end{equation}

Throughout, let $\widehat{\text{ASD}}(\cdot)=\sqrt{\widehat{\text{ASV}}(\cdot)}$ denote the estimated asymptotic standard deviation estimate. Note that the ASV depends on both $f$ and $F$, the density and distribution functions.  There are several options to estimate these, but we choose to use the very flexible Generalized Lambda Distribution (GLD) which, for the FKML parameterization \citep{freimer1988study}, is defined in terms of its quantile function, $Q(p)$,
\begin{equation*}
Q(p)=\lambda_1+\lambda_2^{-1}\left\{\lambda_3^{-1}(p^{\lambda_3} - 1)-\lambda_4^{-1}[(1-p)^{\lambda_4} - 1]\right\} ~,    
\end{equation*}
where $\lambda_1$, $\lambda_2$, $\lambda_3$ and $\lambda_4$ are the location, inverse scale and two shape parameters respectively.
To estimate the GLD parameters we use a recent approach introduced by \cite{dedduwakumara2019efficient} which is computationally efficient making it useful for our simulations that follow.  However, other estimators can also be used.  We then use these parameter estimates with the density and distribution functions for the GLD in R \texttt{gld} package \citep{gld-king}.   

Based on asymptotic normality of the MAD \citep[e.g.][]{falk1997asymptotic}, an asymptotic $(1-\alpha)\%$ confidence interval for $\mad$ is given as 
\begin{equation}
 [L, U]_{\mad} = \Bigg[\widehat{\mad}_X\pm  z_{1-\alpha/2}\;\frac{\wasd{(\MADf,F_{n})}}{\sqrt{n}\,}\Bigg]~,
 \label{eq:CI_MAD}
 \end{equation}
 where the $z_{1-\alpha/2}$ is the $(1-\alpha/2)\times $100 percentile of the standard normal distribution.

When constructing the interval estimator for the squared ratio of $\mad$s, we first derive the confidence interval for the log transformed ratio and then exponentiate to return to the ratio scale. Let $\mathcal{W}(F_1,F_2)=\ln[\R_M(F_1,F_2)]$ then, using the Delta method, it is straightforward to show that $\asv (\mathcal{W},F_1,F_2) \doteq \asv(\R_M,F_1,F_2)/[\R_M(F_1,F_2)]^2$~.  Then a $(1-\alpha)\%$ confidence interval estimator for $R_M$ is given as  
\begin{equation}
 [L, U]_{R_M} = \exp\Bigg[\ln(\widehat{R}_M)\pm  z_{1-\alpha/2}\;\frac{\wasd{(\R_M,F_{n_1},F_{n_2})}}{\widehat{R}_M\sqrt{n_1+n_2}\,}\Bigg]~,
 \label{eq:CI_Rm}
 \end{equation}
where $\widehat{R}_M$ is the squared ratio of MADs estimator and the $\asv$ is in \eqref{eq:ASV_Rm}. 

Finally, a $(1-\alpha)\%$ confidence interval for the difference in MADs is simply
\begin{equation}
 [L, U]_{D_M} = \widehat{D}_M\pm  z_{1-\alpha/2}\;\frac{\wasd{(\D_M,F_{n_1},F_{n_2})}}{\sqrt{n_1+n_2}\,}~,
 \label{eq:CI_Dm}
 \end{equation}
where  $\widehat{D}_M$ is the difference of MADs estimator and the $\asv$ can be found in \eqref{eq:ASV_Dm}. 

\section{Simulations and Examples} \label{sec:SimandEx}
We begin by conducting simulations to assess the coverage properties of the interval estimations for data generated from several distributions.  As pointed out earlier, we have used a new estimator of the GLD parameters provided by \cite{dedduwakumara2019simple} since it exhibits very good performance and is very efficient making it useful for our simulations.  In Appendix A.2, we provide R code for the interval estimators using readily available estimators for the GLD from the \texttt{gld} package \citep{gld-king}.  In that code we have opted for Titterington's method \citep{titterington1985comment} since it to has good performance, albeit is more time consuming.
\subsection{Simulations} \label{sec:Sim}
To investigate the performance of the MAD, squared ratio of $\mad$s and difference of $\mad$s intervals we consider simulated coverage probability and the average confidence interval width as performance measures. We have selected the log normal (LN), exponential (EXP), chi-square ($\chi^2_5$) and Pareto (PAR) distributions with different sample sizes of $n=50, 100, 200, 500, 1000$.  Each simulation consists of 10,000 trials.

\begin{table}[h!tbp]
  \centering
  \caption{Simulated coverage probabilities (and widths in parentheses) for the 95\% confidence interval for the $\mad$ (* denotes median width reported due to excessively large widths for a small number of intervals that skew the mean).}
    \begin{tabular}{ccccc}
    \toprule
     Sample size & \multicolumn{1}{l}{$X\sim$ LN(0,1)} & $X\sim $ EXP(1)  & $X\sim \chi^2_5$ & $X\sim$ PAR(1,7)  \\
    \midrule
    True $\mad=$ & 0.599 & 0.481 & 1.895 & 0.075\\
    50    & 0.938 (1.43) & 0.936 (1.93) & 0.927 (1.25*) & 0.939 (0.34) \\
    100   & 0.940 (0.37) & 0.939 (0.29) & 0.938 (0.91) & 0.939 (0.05) \\
    200   & 0.938 (0.26) & 0.947 (0.20) & 0.942 (0.65) & 0.944 (0.03) \\
    500   & 0.945 (0.16) & 0.948 (0.12) & 0.947 (0.41) & 0.949 (0.02) \\
    1000  & 0.946 (0.12) & 0.951 (0.09) & 0.944 (0.29) & 0.947 (0.01) \\
    \bottomrule
    \end{tabular}%
  \label{tab:CP_MAD}%
\end{table}%

Simulated coverages and widths for the interval estimator of $\mad$s, from \eqref{eq:CI_MAD}, are provided in Table \ref{tab:CP_MAD} for several distributions. The coverage probabilities are all close to the nominal level of 0.95, even for $n=50$ where coverages were approximately in the vicinity of 0.93-0.94.  Coverages become closer to the nominal level as the sample size increases and, as expected the interval widths decrease with increasing sample size.

\begin{table}[h!tbp]
  \centering
  \caption{Simulated coverage probabilities (and widths in parentheses) for the 95\% confidence interval for the squared ratio of $\mad$s ($R_M$) and difference of $\mad$s ($D_M$) (* Median width reported due to excessively large widths for a small number, between 1\% and 2\%, of intervals).}
    \begin{tabular}{cccccc}
    \toprule
    \multicolumn{1}{c}{Sample sizes} &       & \multicolumn{1}{l}{$X\sim$ LN(0,1)} & $X\sim $ EXP(1)  & $X\sim \chi^2_5$ & $X\sim$ PAR(1,7)  \\
    \multicolumn{1}{c}{($n_1$,$n_2$)} & Measure & $Y\sim $ LN(0,1) & $Y\sim$ EXP(1) & $Y\sim \chi^2_2$    & $Y\sim$ PAR(1,3)  \\
    \midrule
    & True $R_M=$    & 1 & 1 & 3.876 & 0.148 \\
          & True $D_M=$  & 0 & 0 & 0.932 & -0.119  \\
    \midrule
    50,50 & $R_M$    & 0.958 (3.71*) & 0.971 (4.03*) & 0.955 (12.14*) & 0.978 (0.91*) \\
          & $D_M$    & 0.967 (2.55) & 0.972 (3.49) & 0.956 (1.54*) & 0.967 (1.17) \\
    \midrule
    100,100 & $R_M$    & 0.949 (2.23) & 0.958 (1.87*) & 0.954 (6.48*) & 0.960 (0.33*) \\
          & $D_M$    & 0.954 (0.52) & 0.958 (0.42) & 0.952 (1.08) & 0.951 (0.16) \\
    \midrule
    200,200 & $R_M$    & 0.953 (1.37) & 0.946 (1.28) & 0.950 (4.51) & 0.952 (0.22) \\
          & $D_M$    & 0.945 (0.37) & 0.950 (0.28) & 0.950 (0.76) & 0.947 (0.10) \\
    \midrule
    200,500 & $R_M$    & 0.946 (1.09) & 0.951 (1.02) & 0.950 (3.47) & 0.952 (0.17) \\
          & $D_M$    & 0.945 (0.31) & 0.951 (0.23) & 0.946 (0.69) & 0.956 (0.07) \\
    \midrule
    500,500 & $R_M$    & 0.946 (0.81) & 0.952 (0.75) & 0.949 (2.69) & 0.950 (0.12) \\
          & $D_M$    & 0.948 (0.23) & 0.953 (0.17) & 0.950 (0.48) & 0.947 (0.06) \\
    \midrule
    500,1000 & $R_M$    & 0.947 (0.69) & 0.952 (0.64) & 0.948 (2.23) & 0.951 (0.10) \\
          & $D_M$    & 0.947 (0.20) & 0.949 (0.15) & 0.949 (0.45) & 0.948 (0.04) \\
    \midrule
    1000,1000 & $R_M$    & 0.947 (0.56) & 0.949 (0.52) & 0.949 (1.87) & 0.950 (0.09) \\
          & $D_M$    & 0.944 (0.16) & 0.950 (0.12) & 0.952 (0.34) & 0.948 (0.04) \\
    \bottomrule
    \end{tabular}%
  \label{tab:CP_RmDm}%
\end{table}%

Simulated coverages for interval estimators of squared ratio of $\mad$s  and  difference of $\mad$s are provided in Table \ref{tab:CP_RmDm} for several distributions. Results show excellent coverages compared to the coverages of F-test \citep[the coverage probabilities for interval estimator of the $F$-test can be found in Table 3 of][]{arachchige2019interval} which are poor due to the violation of underlying normality assumptions). Coverages are very close to the nominal 0.95 for both the squared $\mad$ ratio and difference of $\mad$ for all the selected distributions, including smaller sample sizes. There are some slightly conservative coverages only for $n=50$ and for other sample sizes the coverages become very close.  For smaller sample sizes a very small number of the intervals were very wide (between 1\% and 2\%) so we report the median width instead.  

\subsection{Prostate data example} \label{sec:Ex1}

The prostate data set, which is available in the \texttt{depthTools} package \citep{dt}, is a normalized subset of the \cite{singh2002gene} prostate data set. The data consists of gene expressions for the 100 most variable genes for 25 normal and 25 tumoral samples. 

\begin{figure}[h!t]
 \centering
  \includegraphics[scale = 0.75]{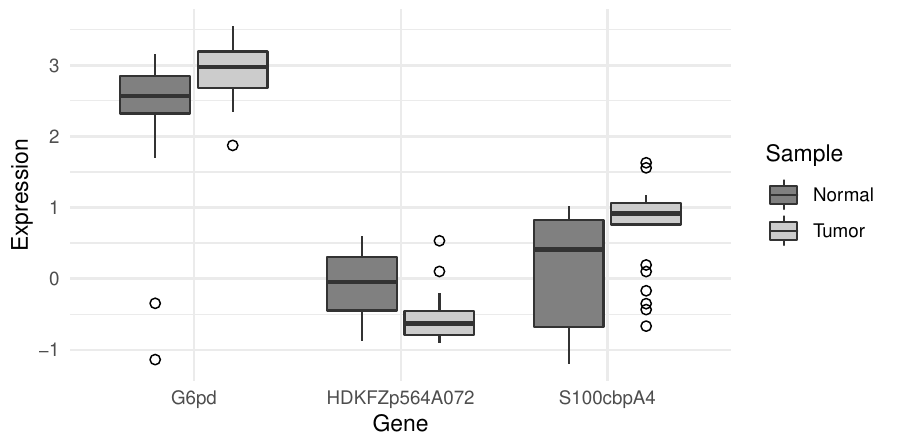} 
  \caption{Box plots of three interesting genes selected from the prostate  data set.}
  \label{fig:gene}
\end{figure}

We selected three genes that are interesting when comparing intervals for ratios of variances and those based on the MADs. These three genes are three of the six that were considered by \cite{arachchige2019interval}. The genes and their abbreviations we consider are Glucose-6-phosphate dehydrogenase (G6pd), HDKFZp564A072 and calcium-binding protein A4 (S100cbpA4). Box plots of the genes are provided in Figure \ref{fig:gene} where we note that, ignoring extreme outliers, the spread for the bulk of the data looks similar for G6pd and very different for HDKFZp564A072 and S100cbpA4.

\begin{table}[htbp]
  \centering
  \caption{95\% asymptotic confidence intervals (CI) for the ratios of variances resulting from the $F$-test ($F$), the ratio of mean absolute deviations (R) \citep{bonett2003confidence}, the ratio of $\mad$s ($R_M$) and difference of $\mad$s ($D_M$) for the three selected genes.}
    \begin{tabular}{cccccc}
    \toprule
    Gene  &       & F     & R     & R\_M  & D\_M \\
    \midrule
    G6pd  & Est.  & 6.496 & 1.723 & 1.000 & 0.000 \\
          & CI    & (2.863, 14.742) &  (0.868, 3.421) & (0.268,  3.734) &  (-0.185, 0.185) \\
    \midrule
    HDKFZp- & Est.  & 1.930 & 1.704 & 5.013 & 0.213 \\
    564A072 & CI    & (0.850, 4.379) &  (1.091, 2.662) &  (1.211, 20.761) &  ( 0.035, 0.391) \\
    \midrule
    S100cbpA4 & Est.  & 1.748 & 1.624 & 8.725 & 0.301 \\
          & CI    & (0.770, 3.968) &  (0.911, 2.894) &  (1.440, 52.856) & (-0.013, 0.615) \\
    \bottomrule
    \end{tabular}%
  \label{tab:prost_CI}%
\end{table}%

 In Table \ref{tab:prost_CI} we provide the point estimate and asymptotic 95\% confidence intervals for the ratio of variances (from the $F$-test assuming underlying normality), the squared ratios of $\mad$s and difference of $\mad$s for the three selected genes. When ignoring two outliers for G6pd the spread looks similar, however the interval for the ratio of variances suggests a large difference in variance between the two.  This is not the case for the MAD intervals where the point and intervals estimates suggest little difference.  For HDKFZp564A072 and S100cbpA4 the intervals tell a different story.  The ratio of variance intervals do not find a significant difference, while the MAD intervals do, or in the case of the difference very close to.  We favor the findings from the MAD due to the obvious difference in spread for the bulk of the data as depicted in the box plots.  This difference in findings is likely due to the group with smaller spread for most data, have extreme outliers that increases the sample variance so that it is similar to the sample variance for the other group. The MADs are not affected by these outliers. \cite{arachchige2019interval} provide similar contrasting results when comparing an asymptotic interval for the ratio of variances and intervals based on the interquantile range.

\section{Summary and discussion} \label{sec:Summary}
The MAD is a  robust estimator of scale exhibiting good robustness properties. We have considered interval estimators for the MAD, ratios of $\mad$s and differences of $\mad$s. Simulation results for the interval estimators showed excellent coverages even for small sample sizes such as $n=50$ for all distributions we considered. Our example reveals that different conclusions can be made by using ratios of $\mad$s and differences of $\mad$s compared to intervals for the ratio of variances which is influenced by outliers.  Future extensions to this work would be to consider intervals for alternatives to the MAD \citep[e.g. see][]{rousseeuw1993alternatives}.

\appendix
\section{Appendix} \label{sec:Appendix}


\subsection{Proof of Theorem \ref{th:PIF_Rm}}\label{app:proof_of_PIF_Rm}
\begin{proof}
 A power series expansion of $\MADf(F_\epsilon)$ can be written as $$\MADf(F_\epsilon)=\MADf(F) + \epsilon \IF(x;\MADf,F)+O(\epsilon^2)~.$$  Let $F_\epsilon=(1-\epsilon)F_1+\epsilon\Delta_{x}$, then we have 
\begin{equation*}
\left[\MADf(F_\epsilon)\right]^2=\MADf^2 (F_1) + 2\epsilon \MADf(F_1)\IF(x;\MADf,F_1)+ O(\epsilon^2)~.
\label{eq:MAD2}
\end{equation*} 

Therefore, the first PIF is
\begin{small}
\begin{align*}
\PIF_1 (x; \R_M,F_1,F_2)=&\lim_{\epsilon \downarrow 0}\Bigg\{\frac{\MADf^2 (F_1) + 2\epsilon \MADf(F_1)\IF(x;\MADf,F_1)+ O(\epsilon^2) -\MADf^2 (F_1)}{\epsilon \MADf^2 (F_2)}\Bigg\}\\
\end{align*}
\end{small}
For the second PIF set $F_\epsilon=(1-\epsilon)F_2+\epsilon\Delta_{x}$. Then 
\begin{align*}
\PIF_2 (x; \R_M,F_1,F_2)=&\lim_{\epsilon \downarrow 0}\Bigg\{\frac{\MADf^2 (F_1)\left[\MADf^2 (F_\epsilon)\right]^{-1} -\MADf^2 (F_1)/\MADf^2 (F_2)}{\epsilon}\Bigg\}\\
=&\lim_{\epsilon \downarrow 0}\Bigg\{\frac{\MADf^2 (F_1)\MADf^2 (F_2) -\MADf^2 (F_1)\MADf^2 (F_\epsilon)}{\epsilon\MADf^2 (F_2)\MADf^2 (F_\epsilon)}\Bigg\}\\
=&\lim_{\epsilon \downarrow 0}\Bigg\{\frac{-2\epsilon\MADf^2 (F_1)\MADf(F_2)\IF(x;\MADf,F_2) + O(\epsilon^2)}{\epsilon\MADf^2 (F_2)\MADf^2 (F_\epsilon)}\Bigg\}
\end{align*}
Recall the $\IF(x;\MADf,F)$ in \eqref{eq:IFmad3} and evaluated at $F_1$ and $F_2$. Finally, the $\PIF_1$ and $\PIF_2$ can be obtained by taking the limit by noting that $\lim_{\epsilon\downarrow 0}[O(\epsilon^2)/\epsilon]=0$.
\end{proof}

\subsection{Bonett and Seier Method }

\cite{bonett2003confidence} suggested constructing distribution-free confidence intervals for median absolute deviation from a target by applying the usual confidence interval for the median described on page no. 137 of \cite{snedecor1980statistical} to the transformed values $|Y_i-h|$.  Note here that the target, $h$, is known and fixed so that this differs from the MAD where the target is the median which needs to be estimated.

\cite{snedecor1980statistical}  describe a simple way to construct the confidence interval for the population median. Let $X_{(1)} \leq X_{(2)} \leq \ldots \leq X_{(n)}$ be the ordered random sample.  Then two of the order statistics become the lower and upper bounds of the confidence interval. The $(1-\alpha)\times 100$\% confidence interval for position of the population median is,
 \begin{equation}
 (n+1)/2 \pm z_{(1-\alpha/2)} \sqrt{n}.  
 \end{equation}

 After rounding the lower limit and upper limit to the nearest integer, two order statistics can be selected as the lower and upper bounds for the confidence interval for the population median. When $X_i=|Y_i-M|$ where $M$ is the true population median and $Y_i$ is the actual value, this confidence interval becomes the interval for the median absolute deviation with the exception that $M$ is assumed fixed and not estimated.

\cite{bonett2003confidence} also suggested constructing distribution-free confidence intervals for a ratio of median absolute deviations from a target by applying the  \cite{price2002distribution} method for transformed values $|Y_{ij}-h|$. The two populations medians were considered as target values. The \cite{price2002distribution} method can be found in detail in Section 2.1 of \cite{arachchige2019interval}.

\begin{table}[htbp]
   \centering
   \caption{Simulated coverage probabilities (and widths in parentheses) for the 95\% confidence interval for the $\mad$ when the target value is the population median and is assumed known.}
     \begin{tabular}{ccccc}
    \toprule
     Sample size & \multicolumn{1}{c}{$X\sim$ LN(0,1)} & \multicolumn{1}{c}{$X\sim$ EXP(1)} & \multicolumn{1}{c}{$X\sim \chi^2_5$} & \multicolumn{1}{c}{$X\sim$ PAR(1,7)} \\
     \midrule
     50    & 0.940(0.31) & 0.937(0.23) & 0.937(1.07) & 0.932(0.04) \\
     100   & 0.944(0.22) & 0.940(0.17) & 0.937(0.79) & 0.946(0.03) \\
    200   & 0.940(0.15) & 0.941(0.12) & 0.939(0.56) & 0.942(0.02) \\
     500   & 0.945(0.10) & 0.946(0.08) & 0.946(0.36) & 0.944(0.01) \\
    1000  & 0.946(0.07) & 0.944(0.05) & 0.946(0.25) & 0.950(0.01) \\
     \bottomrule
     \end{tabular}%
   \label{tab:CI_MADfromTarget1}%
 \end{table}%
From Table \ref{tab:CI_MADfromTarget1}, the  coverage probabilities are close to the nominal coverage (0.95) even for small sample sizes under this assumption of known population median.

\begin{table}[htbp]
  \centering
  \caption{Simulated coverage probabilities (and widths in parentheses) for the 95\% confidence interval for the $\mad$ when the target value is the estimated median.}
    \begin{tabular}{cllll}
    \toprule
    Sample size & \multicolumn{1}{c}{$X\sim$ LN(0,1)} & \multicolumn{1}{c}{$X\sim$ EXP(1)} & \multicolumn{1}{c}{$X\sim \chi^2_5$} & \multicolumn{1}{c}{$X\sim$ PAR(1,7)} \\
    \midrule
    50    & 0.751(0.31) & 0.765(0.24) & 0.900(1.09) & 0.723(0.04) \\
    100   & 0.756(0.22) & 0.773(0.17) & 0.908(0.80) & 0.746(0.03) \\
    200   & 0.755(0.15) & 0.777(0.12) & 0.912(0.56) & 0.737(0.02) \\
    500   & 0.753(0.10) & 0.781(0.08) & 0.912(0.36) & 0.734(0.01) \\
    1000  & 0.754(0.07) & 0.790(0.05) & 0.911(0.25) & 0.740(0.01) \\
    \bottomrule
    \end{tabular}%
  \label{tab:CI_MADfromTarget2}%
\end{table}%

 Table \ref{tab:CI_MADfromTarget2}, shows the simulated coverage probabilities for the confidence interval for the $\mad$ when the target is the estimated median. The coverage probabilities are typically much lower than the nominal coverage of 0.95 even for large sample sizes under this assumption of estimated population median for these selected distributions.  
 
  \begin{table}[htbp]
   \centering
   \caption{Simulated coverage probabilities (and widths in parentheses) for the 95\% confidence interval for the ratio of $\mad$s when the target values are the two population medians.}
    \begin{tabular}{ccccc}
    \toprule
  \multicolumn{1}{c}{Sample sizes} & \multicolumn{1}{l}{$X\sim$ LN(0,1)} & $X\sim $ EXP(1)  & $X\sim \chi^2_5$ & $X\sim$ PAR(1,7)  \\
     \multicolumn{1}{c}{($n_1$,$n_2$)} & $Y\sim $ LN(0,1) & $Y\sim$ EXP(1) & $Y\sim \chi^2_2$    & $Y\sim$ PAR(1,3)  \\
     \midrule
     50,50 & 0.906(0.93) & 0.892(0.89) & 0.900(1.81) & 0.889(0.33) \\
     100,100 & 0.918(0.59) & 0.910(0.57) & 0.912(1.20) & 0.904(0.21) \\
     200,200 & 0.925(0.39) & 0.913(0.38) & 0.918(0.81) & 0.916(0.14) \\
     200,500 & 0.878(0.31) & 0.871(0.30) & 0.859(0.65) & 0.861(0.11) \\
    500,500 & 0.930(0.24) & 0.929(0.24) & 0.927(0.50) & 0.929(0.09) \\
     500,1000 & 0.901(0.20) & 0.896(0.20) & 0.879(0.43) & 0.895(0.07) \\
     1000,1000 & 0.934(0.17) & 0.935(0.16) & 0.929(0.35) & 0.930(0.06) \\
    \bottomrule
     \end{tabular}%
   \label{tab:CI_ratioMADfromTarget}%
 \end{table}%
 
 According to the results in Table \ref{tab:CI_ratioMADfromTarget}, it seems this method works only for moderate to large and equal sample sizes. The coverage probability is bit low for unequal sample sizes.

\subsection{R code for interval estimators}\label{app:proof_of_PIF_Rm}
\begin{verbatim}
# This codes uses the gld R package for estimation of the GLD since it is
# readily available in R.
library(gld)

mad <- function(x) median(abs(x - median(x)))

asv.mad <- function(x, method = "TM"){
  lambda <- fit.fkml(x, method = method)$lambda
  m  <- median(x)
  mad.x <- mad(x)
  fFinv <- dgl(c(m - mad.x, m + mad.x, m), lambda1 = lambda)
  FFinv <- pgl(c(m - mad.x, m + mad.x), lambda1 = lambda)
  A <- fFinv[1] + fFinv[2]
  C <- fFinv[1] - fFinv[2]
  B <- C^2 + 4*C*fFinv[3]*(1 - FFinv[2] - FFinv[1])
  (1/(4 * A^2))*(1 + B/fFinv[3]^2) 
} 

ci.mad <- function(x, y = NULL, gld.est = "TM", 
          two.samp.diff = TRUE, conf.level = 0.95){
  alpha <- 1 - conf.level
  z <- qnorm(1 - alpha/2)
  x <- x[!is.na(x)]
  est <- mad.x <- mad(x)
  n.x <- length(x)
  asv.x <- asv.mad(x, method = gld.est)
  if(is.null(y)){
    ci <- mad.x + c(-z, z)*sqrt(asv.x/n.x)
  } else{
    y <- y[!is.na(y)]
    mad.y <- mad(y)
    n.y <- length(y)
    asv.y <- asv.mad(y, method = gld.est)
    if(two.samp.diff){
      est <- mad.x - mad.y 
      ci <- est + c(-z, z)*sqrt(asv.x/n.x + asv.y/n.y)
    } else{
      est <- (mad.x/mad.y)^2
      log.est <- log(est)
      var.est <- 4 * est * ((1/mad.y^2)*asv.x/n.x + (est/mad.y^2)*asv.y/n.y)
      Var.log.est <- (1 / est^2) * var.est
      ci <- exp(log.est + c(-z, z) * sqrt(Var.log.est))
    }
  }
  list(Estimate = est, conf.int = ci)
}

x <- rlnorm(100)
y <- rlnorm(200, meanlog = 1.2)

ci.mad(x) # single sample
ci.mad(x, y) # two sample difference
ci.mad(x, y, two.samp.diff = FALSE) # two sample squared ratio
\end{verbatim}

\bibliographystyle{authordate4}
\bibliography{ref}

\end{document}